\documentclass[11pt]{amsart}
\usepackage{latexsym}
\usepackage{amssymb,amscd,amsthm,amsfonts,verbatim,amsopn,color,url}
\usepackage[dvips]{epsfig}
\def\FF{{\mathcal F}}
\def\GG{{\mathcal G}}
\def\HH{{\mathcal H}}
\def\II{{\mathcal I}}
\def\LL{{\mathcal L}}

\def\NN{{\mathcal N}}

\def\Ss{{\mathcal S}}
\def\TT{{\mathcal T}}

\newcommand{\dc}{\mathbf{dc}_{\omega}}
\newcommand{\lra}{\longrightarrow}

\newcommand{\mfA}{\mathfrak A}
\newcommand{\nbc}{\mathbf{nbc}_{\omega}}
\newcommand{\nbclex}{\mathbf{nbc}_{lex}}

\newcommand{\pn}{\mathbf{pn}_{\omega}}
\newcommand{\ra}{\rightarrow}

\newcommand{\lf}{\lfloor}
\newcommand{\rf}{\rfloor}

\newcommand{\wti}{\widetilde}

\newtheorem{thm}{Theorem}[section]
\newtheorem{df} [thm]{Definition}

\newtheorem{crl} [thm]{Corollary}
\newtheorem{prop}[thm]{Proposition}

\newtheorem{rem}[thm]{Remark}
\newtheorem{rems}[thm]{Remarks}

\newtheorem{expl}[thm]{Example}
\newtheorem{qst}[thm]{Question}
\newtheorem{qsts}[thm]{Questions}
\numberwithin{equation}{section}

\newenvironment{pf}{\noindent {\bf Proof.}}{\hfill $\Box$\vspace{0.3cm}}
\newenvironment{explrm}{\begin{expl} \rm}{\end{expl}}
\newenvironment{remrm}{\begin{rem} \rm}{\end{rem}}
\newenvironment{remsrm}{\begin{rems} \rm}{\end{rems}}
\newenvironment{dfrm}{\begin{df} \rm}{\end{df}}
\newenvironment{qstrm}{\begin{qst} \rm}{\end{qst}}

\begin{document}

\title[Complexes of trees and nested set complexes]
{Complexes of trees and nested set complexes}

\author{Eva Maria Feichtner}

\address{
Department of Mathematics, ETH Zurich, 8092 Zurich, Switzerland}
\email{feichtne@math.ethz.ch} 

\date{September 2005}
 % general topology of complexes.

\begin{abstract}
  We exhibit an identity of abstract simplicial complexes between the
  well-studied complex of trees $T_n$ and the reduced minimal nested set
  complex of the partition lattice. We conclude that the order complex
  of the partition lattice can be obtained from the complex of trees
  by a sequence of stellar subdivisions. 
  We provide an explicit cohomology basis for the complex of
  trees that emerges naturally from this context.

  Motivated by these results, we review the generalization of
  complexes of trees to complexes of $k$-trees by Hanlon, and we
  propose yet another, in the context of nested set complexes more 
  natural, generalization.
\end{abstract}

%%%%%%%%%%%%%%%%%%%%%%%%%%%%%%%%%%%%%%%%%%%%%%%%%%%%%%%%%%%%%%%%%%%%%%%
%%%%%%%%%%%%%%%%%%%%%%%%%%%%%%%%%%%%%%%%%%%%%%%%%%%%%%%%%%%%%%%%%%%%%%%
%%%%%%%%%%%%%%%%%%%%%%%%%%%%%%%%%%%%%%%%%%%%%%%%%%%%%%%%%%%%%%%%%%%%%%%

\maketitle

%%%%%%%%%%%%%%%%%%%%%%%%%%%%%%%%%%%%%%%%%%%%%%%%%%%%%%%%%%%%%%%%%%%%%%%
%%%%%%%%%%%%%%%%%%%%%%%%%%%%%%%%%%%%%%%%%%%%%%%%%%%%%%%%%%%%%%%%%%%%%%%
%%%%%%%%%%%%%%%%%%%%%%%%%%%%%%%%%%%%%%%%%%%%%%%%%%%%%%%%%%%%%%%%%%%%%%%

\section{Introduction}
\label{sec_intr}

In this article we explore the connection between complexes of trees and 
nested set complexes of specific lattices.

{\em Nested set complexes\/} appear as the combinatorial core in
De~Concini-Procesi wonderful compactifications of arrangement
complements~\cite{DP1}. They record the incidence structure of natural
stratifications and are crucial for descriptions of topological
invariants in combinatorial terms.  Disregarding their geometric
origin, nested set complexes can be defined for any finite
meet-semilattice~\cite{FK}. Interesting connections between seemingly
distant fields have been established when relating the purely
oder-theoretic concept of nested sets to various contexts in geometry.
See~\cite{FY} for a construction linking nested set complexes to toric
geometry, and~\cite{FS} for an appearance of nested set complexes in
tropical geometry.

This paper presents yet another setting where nested set complexes
appear in a meaningful way and, this time, contribute to the toolbox of 
topological combinatorics and combinatorial representation theory.

{\em Complexes of trees\/} $T_n$ are abstract simplicial complexes
with simplices corresponding to combinatorial types of rooted trees on
$n$ labelled leaves.  They made their first appearance in work of
Boardman~\cite{B} in connection with $E_{\infty}$-structures in
homotopy theory. Later, they were studied by Vogtman~\cite{V} from the
point of view of geometric group theory, and by Robinson and
Whitehouse~\cite{RW} from the point of view of representation theory.
In fact, $T_n$ carries a natural action of the symmetric group
$\Sigma_n$ that allows for a lifting to a $\Sigma_{n+1}$-action. For
studying induced representations in homology, Robinson and Whitehouse
determined the homotopy type of $T_n$ to be a wedge of $(n{-}1)!$
spheres of dimension $n{-}3$. Later on, complexes of trees were shown
to be shellable by Trappmann and Ziegler~\cite{TZ} and, in independent
(unpublished) work, by Wachs~\cite{W1}. Recent interest in the
complexes is motivated by the study of spaces of phylogenetic trees
from combinatorial, geometric and statistics point of view~\cite{BHV}.
Complexes of trees appear as links of the origin in natural polyhedral
decompositions of the spaces of phylogenetic trees.

Ardila and Klivans~\cite{AK} recently proved that the complex of trees
$T_n$ can be subdivided by the order complex of the partition lattice
$\Delta(\Pi_n)$. Our result shows that $\Delta(\Pi_n)$, in fact, can
be obtained by a sequence of stellar subdivisions from the complex of
trees. This and other corollaries rely on the specific
properties of nested set complexes that we introduce into the picture.

\bigskip

Our paper is organized as follows: After recalling the definitions of
complexes of trees and of nested set complexes in
Section~\ref{sec_defn}, we establish an isomorphism between the
complex of trees $T_n$ and the reduced minimal nested set complex of
the partition lattice~$\Pi_n$ in Section~\ref{sec_treesAn}. Among
several corollaries, we observe that the isomorphism provides a
$\Sigma_n$-invariant approach for studying tree complexes; their
$\Sigma_n$-representation theory can be retrieved literally for
free.

In Section~\ref{sec_cohbasis} we complement the by now classical
combinatorial correspondence between no broken circuit bases and
decreasing EL-labelled chains for geometric lattices by incorporating
proper maximal nested sets as recently defined by De~Concini and
Procesi~\cite{DP2}. We formulate a cohomology basis for the complex of
trees that emerges naturally from this combinatorial setting.

Mostly due to their rich representation theory, complexes of trees
have been generalized early on to complexes of homeomorphically
irreducible $k$-trees by Hanlon~\cite{H}. We discuss this and another,
in the nested set context more natural, generalization in
Section~\ref{sec_ktrees}.

\bigskip
\noindent
{\em Acknowledgments:} I would like to thank Michelle Wachs and
Federico Ardila for stimulating discussions at the IAS/Park City
Mathematics Institute in July 2004.

%%%%%%%%%%%%%%%%%%%%%%%%%%%%%%%%%%%%%%%%%%%%%%%%%%%%%%%%%%%%%%%%%%%%%%%
%%%%%%%%%%%%%%%%%%%%%%%%%%%%%%%%%%%%%%%%%%%%%%%%%%%%%%%%%%%%%%%%%%%%%%%
%%%%%%%%%%%%%%%%%%%%%%%%%%%%%%%%%%%%%%%%%%%%%%%%%%%%%%%%%%%%%%%%%%%%%%%

\section{Main Characters}
\label{sec_defn}

\subsection{The complex of trees}
\label{ssec_trees}

Let us fix some terminology: A {\em tree} is a cycle-free graph;
vertices of degree~$1$ are called {\em leaves\/} of the tree. A {\em
  rooted tree\/} is a tree with one vertex of degree larger~$1$ marked
as the root of the tree. Vertices other than the leaves and the root
are called {\em internal\/} vertices. We assume that internal vertices
have degree at least~$3$. The root of the tree is thus the only vertex
that can have degree~$2$. Another way of saying this is that we assume
all non-leaves to have {\em outdegree\/} at least~$2$, where the
outdegree of a vertex is the number of adjacent edges that do not lie
on the unique path between the vertex and the root.

We call a rooted tree {\em binary\/} if the vertex degrees are
minimal, i.e., the root has degree~$2$ and the internal vertices have
degree~$3$; in other words, if the outdegree of all non-leaves
is~$2$. Observe that a rooted binary tree on $n$ leaves has exactly
$n{-}2$ internal edges, i.e., edges that are not adjacent to a leave.

The {\em combinatorial type\/} of a rooted tree with labelled leaves
refers to its equivalence class under label- and root-preserving
homeomorphisms of trees as \mbox{$1$-dim}\-en\-sional cell complexes.

Rooted trees on $n$ leaves labelled with
integers $1,\ldots,n$ are in one-to-one correspondence with trees on
$n{+}1$ leaves labelled with integers $0,\ldots,n$ where all internal
vertices have degree~$\geq 3$. The correspondence is obtained by
adding an edge and a leaf labelled~$0$ to the root of the tree.
Though trees on $n{+}1$ labelled leaves seem to be the more natural,
more symmetric objects, rooted trees
on $n$ leaves come in more handy for the description of $T_n$.

\begin{dfrm}
\label{df_Tn}
The {\em complex of trees\/} $T_n$, $n\,{\geq}\,3$, is the abstract
simplicial complex with maximal simplices given by the combinatorial
types of binary rooted trees with~$n$ leaves labelled $1,\ldots,n$,
and lower dimensional simplices obtained by contracting at most
$n{-}3$ internal edges.
\end{dfrm}

The complex of trees $T_n$ is a pure $(n{-}3)$-dimensional simplicial
complex. As we pointed out in the introduction, it is homotopy
equivalent to a wedge of $(n{-}1)!$ spheres of dimension
$n{-}3$~\cite[Thm.\ 1.5]{RW}.

\begin{figure}[ht]
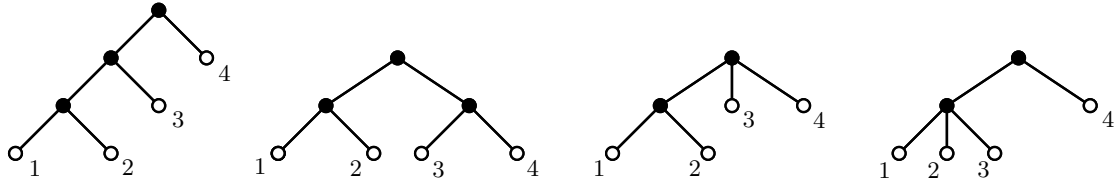

  \begin{picture}(0,0)%
    \includegraphics{T4.pstex}%
  \end{picture}%
  \input{T4.pstex_t}%
  
\caption{Simplices in $T_4$}
\label{fig_T4}
\end{figure}

The complex $T_3$ consists of $3$ points. For $n\,{=}\,4$, there
are~$4$ types of trees, we depict labelled representatives in
Figure~\ref{fig_T4}. Observe that the first two correspond to
$1$-dimensional (maximal) simplices, whereas the other two correspond
to vertices. The last two labelled trees, in fact, are the vertices of
the edge corresponding to the first maximal tree. The third tree is a 
``vertex'' of the second, and of the first. 

%%%%%%%%%%%%%%%%%%%%%%%%%%%%%%%%%%%%%%%%%%%%%%%%%%%%%%%%%%%%%%%%%%%%%%%%%%%%%%

\subsection{Nested set complexes}
\label{ssec_nsc}

We recall here the definition of building sets, nested sets, and
nested set complexes for
finite lattices as proposed in ~\cite{FK}.

We use the standard notation for intervals in a finite lattice~$\LL$,
$[X,Y]\,{:=}\,\{Z\in \LL\,|$ $ X\,{\leq}\,Z\,{\leq}\,Y\}$, for
$X,Y\,{\in}\,\LL$, moreover, $\Ss_{\leq
  X}\,{:=}\,\{Y\,{\in}\,\Ss\,|\,Y\,{\leq}\,X\}$, and accordingly
$\Ss_{< X}$, $\Ss_{\geq X}$, and $\Ss_{> X}$, for
$\Ss\,{\subseteq}\,\LL$ and $X\,{\in}\,\LL$. With $\mathrm{max}\, \Ss$
we denote the set of maximal elements in~$\Ss$ with respect to the
order coming from~$\LL$.

\begin{dfrm} \label{df_building}
Let $\LL$ be a finite lattice. A subset $\GG$ 
in~$\LL_{>\hat 0}$  is called a {\em building set\/} if for any 
$X\,{\in}\,\LL_{>\hat 0}$  and
{\rm max}$\, \GG_{\leq X}=\{G_1,\ldots,G_k\}$ there is an isomorphism
of partially ordered sets
\begin{equation}\label{eq_buildg}
\varphi_X:\,\,\, \prod_{j=1}^k\,\,\, [\hat 0,G_j] \,\, 
                               \stackrel{\cong}{\lra}
                                \,\, [\hat 0,X]
\end{equation}
with $\varphi_X(\hat 0, \ldots, G_j, \ldots, \hat 0)\, = \, G_j$ 
for $j=1,\ldots, k$. 
We call $F_{\GG}(X)\,{:=}\, 
\mathrm{max}\,\GG_{\leq X}$  the {\em set of factors\/} of $X$ in $\GG$.
\end{dfrm}

The full lattice $\LL_{>\hat 0}$ is the simplest example of a building set 
for $\LL$. We will sometimes abuse notation and just write $\LL$ in this case.
Besides this maximal building set, there is always a minimal 
building set $\II$  consisting 
of all elements $X$ in $\LL_{>\hat 0}$ which do not allow for a 
product decomposition of the lower interval~$[\hat 0,X]$, 
the so-called {\em irreducible elements\/} in~$\LL$.

\begin{dfrm} \label{df_nested}
Let $\LL$ be a finite lattice and $\GG$ a building set 
containing the maximal element $\hat 1$ of $\LL$.  A
subset $\Ss$ in $\GG$ is called~{\em nested\/} (or $\GG$-{\em
nested\/} if specification is needed) if, for any set of incomparable
elements $X_1,\dots,X_t$ in $\Ss$ of cardinality at least two, the
join $X_1\vee\dots\vee X_t$ does not belong to $\GG$.  The
$\GG$-nested sets form an abstract simplicial complex, $\wti
\NN(\LL,\GG)$, the {\em nested set complex\/} of $\LL$ with respect 
to~$\GG$. Topologically, the nested set complex is a cone with 
apex $\hat 1$; its base $\NN(\LL,\GG)$ is called the {\em reduced 
nested set complex\/} of $\LL$ with respect to~$\GG$. 
\end{dfrm}

We will mostly be concerned with reduced nested set complexes due to
their more interesting topology. If the underlying lattice is clear
from the context, we will write $\NN(\GG)$ for $\NN(\LL,\GG)$.

Nested set complexes can be defined analogously for building sets not 
containing~$\hat 1$, and, even more generally, for meet semi-lattices. For a 
definition in the full generality, see~\cite[Section~2]{FK}.

For the maximal building set of a lattice $\LL$, subsets are nested if
and only if they are linearly ordered in $\LL$. Hence, the reduced
nested set complex $\NN(\LL,\LL)$ coincides with the order complex of
$\LL$, more precisely, with the order complex of the proper part,
$\LL\,{\setminus}\,\{\hat 0,\hat 1\}$, of $\LL$, which we denote
by~$\Delta(\LL)$ using customary notation.

If $\LL$ is an atomic lattice, the nested set complexes can be
realized as simplicial fans, see~\cite{FY}, and for building sets
$\GG\,{\subseteq}\,\HH$ in $\LL$, the nested set complex $\wti
\NN(\LL,\HH)$ can be obtained from $\wti \NN(\LL,\GG)$ by a sequence
of stellar subdivisions~\cite[Thm\ 4.2]{FM}. 
In particular, any reduced nested 
set complex $\NN(\LL,\GG)$ is
obtained by a sequence of stellar subdivisions from the minimal
reduced nested set complex $\NN(\LL,\II)$, and can be further
subdivided by stellar subdivisions so as to obtain the maximal nested
set complex $\Delta(\LL)$.

\begin{explrm}
\label{ex_NPin}
Let $\Pi_n$ denote the lattice of set partitions of
$[n]\,{:=}\,\{1,\ldots,n\}$ partially ordered by reversed refinement.  As
explained above, the reduced maximal nested set complex
$\NN(\Pi_n,\Pi_n)$ is the order complex $\Delta(\Pi_n)$. Irreducible
elements in $\Pi_n$ are the partitions with exactly one non-singleton
block. They can be identified with subsets of~$[n]$ of cardinality at
least~$2$. Nested sets for the minimal building set $\II$ are
collections of such subsets of~$[n]$ such that any two either contain
one another or are disjoint. For $n\,{=}\,3$, the reduced minimal
nested set complex consists of $3$ isolated points; for $n\,{=}\,4$,
it equals the Petersen graph.
\end{explrm}

%%%%%%%%%%%%%%%%%%%%%%%%%%%%%%%%%%%%%%%%%%%%%%%%%%%%%%%%%%%%%%%%%%%%%%%
%%%%%%%%%%%%%%%%%%%%%%%%%%%%%%%%%%%%%%%%%%%%%%%%%%%%%%%%%%%%%%%%%%%%%%%
%%%%%%%%%%%%%%%%%%%%%%%%%%%%%%%%%%%%%%%%%%%%%%%%%%%%%%%%%%%%%%%%%%%%%%%

\section{Subdividing the complex of trees}
\label{sec_treesAn}

We now state the core fact of our note.

\begin{thm}
\label{thm_TnminN}
The complex of trees $T_n$ and the reduced minimal nested set complex of the 
partition lattice $\NN(\Pi_n,\II)$ coincide as abstract simplicial 
complexes.
\end{thm}

\begin{pf}
We exhibit a bijection between simplices in $T_n$ and nested sets in
the reduced minimal nested set complex $\NN(\Pi_n,\II)$ of the
partition lattice $\Pi_n$.

Let $T$ be a tree in $T_n$ with inner vertices $t_1, \ldots, t_k$. 
We denote the set of leaves in $T$ below an inner vertex $t$ by $\ell(t)$.
We associate a nested set $\Ss(T)$ in $\NN(\Pi_n,\II)$ to $T$ by defining
\[
     \Ss(T) \, :=\, \{ \ell(t_i)\,|\, i=1,\ldots,k \}\, .      
\] 

Conversely, let $\Ss\,{=}\,\{S_1,\ldots, S_k\}$ be a (reduced) nested
set in $\Pi_n$ with respect to $\II$. We define a rooted tree $\wti
T(\Ss)$ on the vertex set $\Ss\,{\cup}\, \{R\}$, where $R$ will be the
root of the tree. Cover relations are defined by setting $S\,{>}\,T$ if and
only if $T\,{\in}\,{\rm max}\Ss_{<S}$, and requiring the root $R$ to cover 
any element in ${\rm max}\Ss$. Observe that we allow the ``root'' of $\wti
T(\Ss)$ to have degree~$1$ in this intermediate stage of the
construction.

To obtain a
tree with $n$ leaves, i.e., a simplex in $T_n$, we need to ``grow
leaves'' on $\wti T(\Ss)$. To do so, expand any leaf $S$ of $\wti
T(\Ss)$ into $|S|$ many leaves labeled with the elements in $S$. For
every internal vertex $T$ of $\wti T(\Ss)$ add leaves labeled with the
elements of $T\,{\setminus}\,\bigcup_{S\in \Ss, S<T} S$, analogously, add
leaves labeled with $[n]\,{\setminus}\,\bigcup \Ss$ to the root $R$. 
Denote the resulting tree by $T(\Ss)$.

The two maps are inverse to each other, hence we have an
(order-preserving) bijection between the faces of $T_n$ and the reduced
minimal nested set complex $\NN(\Pi_n,\II)$ of the partition lattice. 
We see that $T_n$ and $\NN(\Pi_n,\II)$, in fact, are identical as abstract
simplicial complexes.
\end{pf}

\begin{crl}
\label{crl_subd}
The order complex of the partition lattice $\Delta(\Pi_n)$ 
can be obtained from the complex of trees $T_n$ by a sequence 
of stellar subdivisions. 
\end{crl}

\begin{pf}
Referring to \cite[Thm.\ 4.2]{FM}, we see that the order complex
$\Delta (\Pi_n)$, i.e., the reduced maximal nested set complex
$\NN(\Pi_n,\Pi_n)$ of $\Pi_n$, can be obtained from the complex of
trees $T_n\,{=}\,\NN(\Pi_n,\II)$ by a sequence of stellar
subdivisions.

Explicitly, the subdivision is given by the choice of a linear 
extension order on $(\Pi_n\setminus \II)^{\rm op}$, and by performing 
stellar subdivisions in  simplices $F_{\II}(X)$ for poset elements $X$ 
in $(\Pi_n\setminus \II)^{\rm op}$ along the given linear order.  
\end{pf}

Both the complex of trees and the partition lattice carry a natural
action of the symmetric group $\Sigma_n$, which induces a
$\Sigma_n$-action on the respective homology, each concentrated in top
dimension (compare~\cite{RW} for the complex of trees, and \cite{S}
for the partition lattice). We recover the well-known isomorphism of
$\Sigma_n$-modules in

\begin{crl}
\label{crl_Snaction}
The top degree homology groups of \,$T_n$ and of \,$\Pi_n$ are isomorphic as 
$\Sigma_n$-modules:
\[ 
            \wti H_{n-3}(T_n)\, \, \cong_{\Sigma_n} \,\,  \wti H_{n-3}(\Pi_n)\, .
\]
\end{crl}

\begin{pf}
  The maps defined in the proof of Theorem~\ref{thm_TnminN} are
  $\Sigma_n$-equivariant, hence they induce a $\Sigma_n$-module
  isomorphism in homology. Moreover, the subdivision referred to in
  the proof of Corollary~\ref{crl_subd} is $\Sigma_n$-equivariant,
  where the action on $\Delta(\Pi_n)$ comes from the action of
  $\Sigma_n$ on $\Pi_n$.  Observe that not any individual subdivision
  step is equivariant, but simultaneous subdivisions for a fixed
  partition type are.
\end{pf}

\begin{remsrm}
\noindent
{\bf (1)} In recent work, Ardila and Klivans showed that the complex
of trees $T_n$ coincides as a simplicial complex with the matroid
stratification of the Bergman complex of the graphical matroid
$M(K_n)$~\cite[Sect.\ 3]{AK}. This implies that $\Delta(\LL(M(K_n)))\,
{=}\,\Delta (\Pi_n)$ subdivides the complex of trees, in
fact, the order complex coincides with the weight stratification of
the Bergman complex. Our result gives an explicit sequence of stellar
subdivisions relating one complex to the other. Independently,
Ardila~\cite{A} has described a sequence of subdivisions
connecting~$T_n$ and $\Delta(\Pi_n)$ which remains to be compared to
the one presented here.

\noindent
{\bf (2)} For an arbitrary matroid $M$, the matroid stratification of
the Bergman complex $B(M)$ can be subdivided so as to realize the
order complex of the lattice of flats~$\LL(M)$~\cite[Thm.\ 1]{AK}.
This stratification, in fact, is the finest in a family of
subdivisions having the combinatorics of reduced nested set complexes
$\NN(\LL(M), \GG)$, the coarsest of which, corresponding to
$\NN(\LL(M),\II)$, still subdivides the matroid stratification of
$B(M)$~\cite{FS}.

Our theorem implies that for the particular case of the graphic matroid 
$M(K_n)$ the coarsest nested set stratification coincides with the 
matroid stratification of the Bergman complex. This is not true in 
general; for an example as well as a characterization of when the 
stratifications coincide, see~\cite[Ex.\ 1.2, Thm.\ 5.3]{FS}.

\noindent
{\bf (3)}  The subdivision of $T_n$ described in Corollary~\ref{crl_subd}
can not be connected to the barycentric subdivision of the complex of trees,
${\rm bsd}(T_n)$, by further stellar subdivisions. In earlier work, M.\ Wachs
related cells in $T_n$ to simplices in the barycentric 
subdivision of $T_n$ in a way that suggested such a connection~\cite{W2}. 

As an example, let us observe that the subdivision $\Delta (\Pi_5)$ of
$T_5$ is not refined by ${\rm bsd}(T_5)$:

Consider the triangle $\TT\,{=}\,\{23,45,145\}$ in
$\NN(\Pi_5,\II)\,{=}\,T_5$. The sequence of stellar subdivisions found
in Corollary~\ref{crl_subd} has to be non-increasing on non-building
set elements, hence, in the course of subdividing $T_5$, we first do a
stellar subdivision in the edge $(23,145)$, and later in the edge
$(23,45)$.  The resulting subdivision of $\TT$ is depicted in
Figure~\ref{fig_triangle}.  The edge $(23|45\,{<}\,23|145)$ in
$\Delta(T_5)$ is not refined by the barycentric subdivision of $\TT$.

\begin{figure}[ht]
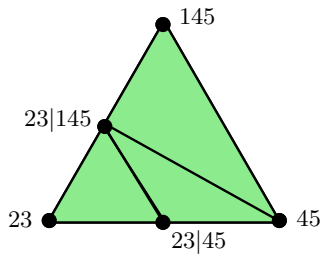

  \begin{picture}(0,0)%
    \includegraphics{triangle.pstex}%
  \end{picture}%
  \input{triangle.pstex_t}%
  
\caption{Subdivision of $\TT\,{=}\,\{23,45,145\}$ in $\Delta(\Pi_5)$}
\label{fig_triangle}
\end{figure}

\end{remsrm}

%%%%%%%%%%%%%%%%%%%%%%%%%%%%%%%%%%%%%%%%%%%%%%%%%%%%%%%%%%%%%%%%%%%%%%%
%%%%%%%%%%   Cohomology basis for $T_n$   %%%%%%%%%%%%%%%%%%%%%%%%%%%%%
%%%%%%%%%%%%%%%%%%%%%%%%%%%%%%%%%%%%%%%%%%%%%%%%%%%%%%%%%%%%%%%%%%%%%%%

\section{A cohomology basis for the complex of trees}
\label{sec_cohbasis}

\noindent 
We start this section by recalling some standard constructions associated 
with geometric lattices. For standard matroid terminology used throughout, 
we refer to the book of J.~Oxley~\cite{Ox}. 

Let $\LL$ be a geometric lattice of rank $r$ with a fixed linear order
$\omega$ on its set of atoms $\mfA(\LL)$. A {\em broken circuit\/} in
$\LL$ is a subset $D$ of $\mfA(\LL)$ that can be written as a circuit
in $\LL$ with its minimal element removed, $D\,{=}\,C\,{\setminus}\min
C$, $C$ a circuit in~$\LL$.  The maximal independent sets in $\LL$ that 
do not contain any broken circuits are known as the {\em no broken
  circuit bases\/} of $\LL$, denoted by~$\nbc(\LL)$.
The no broken circuit bases play an important role for 
understanding the (co)homology of $\Delta(\LL)$.
 
There is a standard labelling $\lambda$ of cover relations $X{>}Y$ in
$\LL$. For $X\,{\in}\,\LL$, let $\lf X \rf\,{:=}\,\{A\in
\mfA(\LL)\,|\,A\leq X \}$, the set of atoms in $\LL$ below $X$, and
define
\[
   \lambda(X>Y)\,\, : = \,\, {\rm min} (\lf X \rf \setminus \lf Y \rf)\, .
\]

This labelling in fact is an $EL$-labelling of $\LL$ in the sense
of~\cite{BWa1}, thus, by ordering maximal chains lexicographically, it
induces a shelling of the order complex $\Delta(\LL)$. Denote by
$\dc(\LL)$ the set of maximal chains in $\LL$ with (strictly)
decreasing label sequence:
\[
    \dc(\LL)\,:=\,\{\hat 0<c_1<\ldots<c_{r-1}<\hat 1\,|\,
    \lambda(c_1>\hat0)>\lambda(c_2>c_1)>\ldots >\lambda(\hat 1>c_{r-1})\}\, .
\]
The characteristic cohomology classes $[c^*]$ for 
$c\,{\in}\,\dc(\LL)$, i.e., classes represented by 
cochains that evaluate to~$1$ on~$c$ and to~$0$ on any other 
top dimensional simplex of~$\Delta(\LL)$, form a basis of the only 
non-zero reduced cohomology group of the order complex, 
$\wti H^{r-1}(\Delta(\LL))$. 

We add the notion of {\em proper\/} maximal nested sets to the
standard notions for geometric lattices with fixed atom order that we
listed so far.  The concept has appeared in recent work of De~Concini
and Procesi~\cite{DP2}. Define a map $\phi:\, \II \,{\rightarrow}\,
\mfA(\LL)$ by setting $\phi(S):={\rm min} \lf S\rf$ for
$S\,{\in}\,\II$.  A maximal nested set $\Ss$ in the (non-reduced)
nested set complex $\wti\NN(\LL,\II)$ is called {\em proper\/} if the
set $\{\phi(S)\,|\, S\in \Ss \,\}$ is a basis of $\LL$. Denote the set
of proper maximal nested sets in~$\LL$ by $\pn(\LL)$.

We define maps connecting $\nbc(\LL)$, $\dc(\LL)$, and $\pn(\LL)$ for 
a given geometric lattice~$\LL$. In the following proposition we will see 
that these maps provide bijective correspondences between the respective 
sets.

To begin with, define $\Psi: \nbc(\LL) \rightarrow \dc(\LL)$ by 
\begin{equation} \label{eq_defPsi}
\Psi(a_1,\ldots,a_r)\,\,=\,\, 
  (\hat 0 < a_r < a_r \vee a_{r-1} < \ldots 
      < a_r \vee a_{r-1} \vee \ldots \vee a_1=\hat 1) \, , 
\end{equation}
where the $a_1,\ldots, a_r$ are assumed to be in ascending order with 
respect to~$\omega$. 

Next, define $\Theta: \dc(\LL) \rightarrow \pn(\LL)$ by
\begin{equation} \label{eq_defTheta}
\Theta(\hat 0<c_1< \ldots < c_{r-1}< \hat 1) \,\, = \,\, 
 F(c_1) \, \cup\, F(c_2)\, \cup\, \ldots  \, \cup\, F(c_{r-1})
        \, \cup\, F(\hat 1)\, ,  
\end{equation} 
for a chain $c:\, \hat 0<c_1< \ldots < c_{r-1}< \hat 1$ in $\LL$ with 
decreasing label sequence, where $F(c_i)$ denotes the set of factors of~$c_i$
with respect to the minimal building set~$\II$ in~$\LL$.

Finally, define $\Phi: \pn(\LL) \rightarrow \nbc(\LL)$ by
\begin{equation} \label{eq_defPhi}
  \Phi(\Ss)\, \, = \, \, \{\phi(S)\,|\, S\in \Ss\}\, ,      
\end{equation}  
for $\Ss\,{\in}\, \pn(\LL)$.

\begin{explrm} \label{ex_bijcor}
Let us consider the partition lattice $\Pi_5$ with lexicographic order 
$lex$
on its set of atoms $ij$, $1\,{\leq}\,i\,{<}\,j\,{\leq}\,5$. 
For $b\,{=}\,\{12,14,23,45\}\,{\in}\,\nbclex(\Pi_5)$ we have 
\[
 \Psi(b)\, \, = \, \, (0<45<23{|}45<23{|}145<\hat1)\, .    
\]
Observe that the chain is constructed by taking consecutive joins of elements
in the opposite of the lexicographic order.

Going further using $\Theta$ we obtain the following proper nested set,
\begin{eqnarray*}
    \Theta(0<45<23{|}45<23{|}145<\hat1) & = &
    \{45\}\, \cup\, \{23,45\}\, \cup\,  \{23,145\}\, \cup\,  \{12345\} \\
& = &
    \{23,45,145,12345\}\, .
\end{eqnarray*}

Applying $\Phi$ we retrieve the no broken circuit basis we started with:
\[
 \Phi(\{23,45,145,12345\})\, \, = \,\, \{12,23,14,45\}\, .       
\]
  
\end{explrm}

\begin{prop} \label{prop_bijcor}
For a geometric lattice $\LL$ with a given linear order $\omega$ on 
its atoms the maps $\Psi$, $\Theta$, and $\Phi$ defined above give
bijective correspondences between (1) the no broken circuit bases 
$\nbc(\LL)$ of $\LL$, (2) the maximal chains in $\LL$ with decreasing 
label sequence, $\dc(\LL)$, and 
(3) the proper maximal nested sets $\pn(\LL)$ 
in $\wti \NN(\LL,\II)$, respectively. 
\end{prop}

\begin{pf}
  The map $\Psi: \nbc(\LL) \rightarrow \dc(\LL)$ is well known in the
  theory of geometric lattices. It is the standard bijection relating
  no broken circuit bases to cohomology generators of the lattice,
  compare \cite[Sect.\ 7.6]{Bj2} for details. In the previously cited
  work of De~Concini and Procesi the composition of maps $\eta:=
  \Theta \circ \Psi: \nbc(\LL) \rightarrow \pn(\LL)$ is shown to be a
  bijection with inverse $\Phi: \pn(\LL) \rightarrow
  \nbc(\LL)$~\cite[Thm.\ 2.2]{DP2}.  This implies that $\Theta:
  \dc(\LL) \rightarrow \pn(\LL)$ is bijective as well, which completes
  the proof of our claim.
\end{pf}

The aim of the next proposition is to trace the support simplices for
the cohomology bases $\{\,[c^*]\,|\, c\,{\in}\,\dc(\Pi_n)\,\}$ of
$\Delta(\Pi_n)$ through the inverse stellar subdivisions linking
$\Delta(\Pi_n)$ to the complex of trees $T_n\,{=}\,\NN(\Pi_n,\II)$.
For the moment we can stay with the full generality of geometric
lattices and study support simplices for maximal simplices of
$\Delta(\LL)$ in the minimal reduced nested set complex
$\NN(\LL,\II)$.

\begin{prop} \label{prop_spts}
  Let $\LL$ be a geometric lattice, $c:\,
  c_1\,{<}\,\ldots\,{<}\,c_{r-1}$ a maximal simplex in $\Delta(\LL)$.
  The maximal simplex in $\NN(\LL,\II)$ supporting $c$ is
  given by the union of sets of factors
\[
    F(c_1)\, \cup \,F(c_2)\, \cup \, \ldots \, \cup \, F(c_{r-1})\, .
\]
\end{prop}

\begin{proof}
There is a sequence of building sets
\[
        \LL\, = \, \GG_{1} \, \supseteq \, \GG_2 \, 
                   \supseteq \, \ldots \,  \supseteq \, \GG_t \, = \, \II\, ,
\]
connecting $\LL$ and $\II$ which is obtained by removing elements of 
$\LL\,{\setminus}\,\II$ from $\LL$ in a non-decreasing order:
$\GG_i\,{\setminus}\,\GG_{i+1}\,{=}\, \{G_i\}$ with $G_i$ minimal in 
$\GG_i\,{\setminus}\,\II$ for $i\,{=}\,1,\ldots,t{-}1$. The corresponding 
nested set complexes are linked by inverse stellar subdivisions: 
\[
      \NN(\GG_i) \, \, = \, \, 
        {\rm st}(\NN(\GG_{i+1}), V(F_{\II}(G_i)))\, ,
      \qquad \mbox{ for }\,i=1,\ldots t{-}1\, ,
\]
where $V(F_{\II}(G_i))$ denotes the simplex in  $\NN(\GG_{i+1})$ spanned
by the factors of $G_i$ with respect to the minimal building set~$\II$.

We trace what happens to the support simplex of $c$ along the sequence of 
inverse stellar subdivisions connecting $\Delta(\LL)$ with $\NN(\II)$.
The support simplex of $c$ remains unchanged in step~$i$ unless $G_i$ 
coincides with a (reducible) chain element $c_j$ (the irreducible chain 
elements can be replaced any time by their ``factors'': 
$F(c_k)\,{=}\,\{c_k\}$  for $c_k\,{\in}\,\II$).

We can assume that the support simplex of $c$ in $\NN(\GG_i)$ is of the form
\[
     \Ss\, = \, F(c_1)\, \cup\, \ldots \,  \cup \, F(c_{j-1})\, 
        \cup \, \{c_j\}\, \cup\, \ldots \, \cup \, \{c_{r-1}\} \, ,
\] 
and we aim to show that  the support simplex of $c$ in $\NN(\GG_{i+1})$ 
is given by
\[
     \TT\, = \, F(c_1)\, \cup\, \ldots \,  \cup \, F(c_{j-1})
                                         \,\cup \, F(c_j)\, 
        \cup \, \{c_{j+1}\}\, \cup\, \ldots \, \cup \, \{c_{r-1}\} \, .
\]
Recall that the respective face posets of the nested set complexes are 
connected by a combinatorial blowup
\begin{equation} \label{eq_combblowup}
   \FF(\NN(\GG_i)) \, \, = \, \, {\rm Bl}_{F(G_i)}
                ( \FF(\NN(\GG_{i+1})))\, .    
\end{equation}
See~\cite[3.1.]{FK} for the concept of a combinatorial blowup in meet
semi-lattices. Hence, the support simplex $\Ss$ of~$c$ in $\NN(\GG_i)$
is of the form $\Ss\,{=}\,\Ss_0\,{\cup}\,\{c_j\}$ with
$\Ss_0\,{\in}\,\NN(\GG_{i+1})$ (it is an element in the ``copy'' of
the lower ideal of elements in $\FF(\NN(\GG_{i+1}))$ having joins with
$F(G_i)$).  Due to~(\ref{eq_combblowup}) we know that
$\Ss_0\,{\not\supseteq}\,F(G_i)$ and
$\Ss_0\,{\cup}\,F(G_i)\,{\in}\,\NN(\GG_{i+1})$, which in fact is the
new support simplex of $c$. Let us mention in passing that, since we
are talking about maximal simplices, $\Ss_0$ contains $F(G_i)$ up to
exactly one element $X_i\,{\in}\,\LL$.

Since $c_j$ is not contained in any of the $F(c_i)$,
$i\,{=}\,1,\ldots,j{-}1$, we have $\Ss_0\,{=}\, F(c_1)\, {\cup}\, $ $
\ldots \, {\cup} \, F(c_{j-1}) \,{\cup} \, \{c_{j+1}\}\, {\cup}\,
\ldots \, {\cup} \, \{c_{r-1}\}$, and we find that $
\TT\,{=}\,F(c_1)\, {\cup}\, \ldots \, {\cup} \, F(c_{j-1}) \,{\cup} \,
F(c_j)\, $ $ {\cup} \, \{c_{j+1}\}\, {\cup}\, \ldots \, {\cup} \,
\{c_{r-1}\}$ as claimed.
\end{proof}

\begin{figure}[ht]
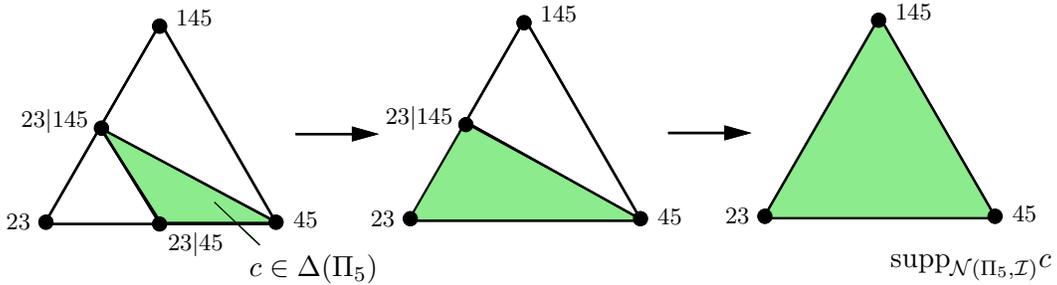

  \begin{picture}(0,0)%
    \includegraphics{subd.pstex}%
  \end{picture}%
  \input{subd.pstex_t}%
  
\caption{Support simplices of $c$}  
\label{fig_subd}
\end{figure}

\begin{explrm} \label{ex_subd}
Let us again consider the partition lattice $\Pi_5$. The support
simplex of $c:\, 45
\,{<}\,23{|}45\,{<}\,23{|}145$ in $\NN(\Pi_5,\II)$ is
$\{23,45,145\}$. 
We depict in Figure~\ref{fig_subd} how the support
simplex of $c$ changes in the sequence of inverse stellar subdivisions
from $\Delta(\Pi_5)$ to $\NN(\Pi_5,\II)$. 
\end{explrm}

We now combine our findings to provide an explicit cohomology basis for 
the complex of trees $T_n$. 

\begin{figure}[ht]
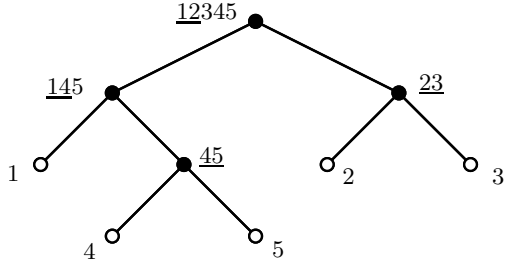

  \begin{picture}(0,0)%
    \includegraphics{admtree.pstex}%
  \end{picture}%
  \input{admtree.pstex_t}%
  
\caption{An admissible tree in $T_5$.}  
\label{fig_admtree}
\end{figure}

We call a binary rooted tree $T$ with  $n$ leaves labelled $1,\ldots,n$ 
{\em admissible\/}, if, when recording the 2nd smallest
label on the sets of leaves below any of the $n{-}1$ non-leaves of $T$, 
we find each of the labels $2,\ldots,n$ exactly once. For an example of an 
admissible tree in $T_5$ see Figure~\ref{fig_admtree}.

\begin{prop}\label{prop_cohbsTn}
\noindent
The characteristic cohomology classes associated with admissible trees 
in $T_n$,
\[
\{\, [T^*]\in \widetilde H^{n-3}(T_n)\,|\, T \mbox{ admissible in }\, 
T_n\,\}\, ,
\]
form a basis for the (reduced) cohomology of the complex of trees $T_n$. 
\end{prop}

\begin{pf}
We set out from the linear basis for $\wti H^{n-3}(\Delta(\Pi_n))$ provided
by characteristic cohomology classes associated with the decreasing
chains $\dc(\Pi_n)$ in $\Pi_n$. Combining Proposition~\ref{prop_spts}
with the definition of the bijection $\Theta:\, \dc(\Pi_n) \ra
\pn(\Pi_n)$ in~(\ref{eq_defTheta}) we find that the characteristic
cohomology classes associated with (reduced) proper maximal nested
sets $\pn(\Pi_n)$ provide a linear basis for
$\wti H^{n-3}(\NN(\Pi_n,\II))$. We tacitly make use of the bijection
between maximal simplices in $\wti \NN(\LL,\II)$ and $\NN(\LL,\II)$ given by
removing the maximal element $\hat 1$ of~$\LL$.

To describe support simplices explicitly, recall that proper maximal
nested sets are inverse images of no broken circuit bases under
$\Phi:\, \pn(\Pi_n) \ra \nbc(\Pi_n)$ as defined in~(\ref{eq_defPhi}).
The no broken circuit bases of $\Pi_n$ with respect to the
lexicographic order on atoms, i.e., on pairs $(i,j)$,
$1\,{\leq}\,i\,{<}\,j\,{\leq}\,n$, are $(n{-}1)$-element subsets of
the form
\[
      (1,2), \, (i_2,3),\, \ldots,\, (i_{n-1},n)
\]
with $1\,{\leq}\,i_j\,{\leq}\,j$ for $j\,{=}\,2,\ldots,n{-}1$. Inverse
images under $\Phi$ are maximal nested sets $\Ss\,{\in}\,\wti
\NN(\Pi_n,\II)$ such that $\{\phi(S)\,|\, S\,{\in}\, \Ss\}$ gives
collections of pairs with each integer from $2$ to $n$ occurring
exactly once in the second coordinate.  Applying the isomorphism
between $\NN(\Pi_n,\II)$ and $T_n$ from Theorem~\ref{thm_TnminN} shows
that the characteristic cohomology classes on admissible trees in
$T_n$ indeed form a basis of $\wti H^{n-3}(T_n)$.
\end{pf}

\begin{remrm}
  Our basis of admissible trees differs from the one presented
  in~\cite[Cor.~5]{TZ} as a consequence of their shelling argument
  for complexes of trees.
\end{remrm}

%%%%%%%%%%%%%%%%%%%%%%%%%%%%%%%%%%%%%%%%%%%%%%%%%%%%%%%%%%%%%%%%%%%%%%%
%%%%%%%%%%%%%%%%%%       $k$-trees     %%%%%%%%%%%%%%%%%%%%%%%%%%%%%%%%
%%%%%%%%%%%%%%%%%%%%%%%%%%%%%%%%%%%%%%%%%%%%%%%%%%%%%%%%%%%%%%%%%%%%%%%

\section{Complexes of $k$-trees and other generalizations}
\label{sec_ktrees}

The intriguing representation theory of complexes of trees $T_n$~\cite{RW} 
has given rise to a generalization to complexes of $k$-trees~\cite{H}.

\begin{dfrm}
\label{df_Tn^k}
The complex of $k$-trees $T_n^{(k)}$, $n\,{\geq}\,1$, $k\,{\geq}\,1$,
is the abstract simplicial complex with faces corresponding to
combinatorial types of rooted trees with~$(n{-}1)k{+}1$ leaves
labelled $1,\ldots,(n{-}1)k{+}1$, with all outdegrees at least $k{+}1$ and
congruent to $1$~mod~$k$, and at least one internal edge. The partial
order among the rooted trees is given by contraction of internal edges.
\end{dfrm}

Alternatively, we could define $T_n^{(k)}$ as the simplicial complex
with faces corresponding to (non-rooted) trees with $(n{-}1)k{+}2$
labelled leaves, all degrees of non-leaves at least $k{+}2$ and
congruent to $2$ mod $k$, and at least one internal edge. Again, the
order relation is given by contracting internal edges.
Observe that for $k\,{=}\,1$ we recover the complex of trees $T_n$.
The face poset of our complex $T_n^{(k)}$ is the poset 
$\LL_{n-1}^{(k)}$ of Hanlon in~\cite{H}.

The complexes $T_n^{(k)}$ are pure simplicial complexes of dimension
$n{-}3$. They were shown to be Cohen-Macaulay by Hanlon~\cite[Thm.\
2.3]{H}; later a shellability result was obtained by Trappmann and
Ziegler~\cite{TZ} and, independently, by Wachs~\cite{W1}.

The complexes $T_n^{(k)}$ carry a natural $\Sigma_N$-action for
$N\,{=}\,(n{-}1)k{+}1$ by permutation of leaves, which induces a
$\Sigma_N$-action on top degree homology. It follows from work of
Hanlon~\cite[Thm.\ 1.1]{H} and Hanlon and Wachs~\cite[Thm.\ 3.11 and
4.13]{HW} that $\wti H_{n-3}(T_n^{(k)})$ is isomorphic as an
$\Sigma_N$-module to $\wti H_{n-3}(\Pi_N^{(k)})$, where $\Pi_N^{(k)}$
is the subposet of $\Pi_N$ consisting of all partitions with block
sizes congruent $1$~mod~$k$. The poset $\Pi_{N}^{(k)}$ had been
studied before on its own right: it was shown to be Cohen-Macaulay by
Bj\"orner~\cite{Bj1}, its homology and $\Sigma_N$-representation
theory was studied by Calderbank, Hanlon and Robinson~\cite{CHR}.

In fact, both $\Sigma_N$-modules $\wti H_{n-3}(T_n^{(k)})$ and $\wti
H_{n-3}(\Pi_N^{(k)})$ are isomorphic to the $1^N$ homogeneous piece of
the free Lie $k$-algebra constructed in~\cite{HW}.

This certainly provides enough evidence to look for a topological
explanation of the isomorphism of $\Sigma_N$-modules:

\begin{qstrm}
\label{qst_ktrees}
Is the complex of $k$-trees $T_n^{(k)}$ related to the order complex
of $\Pi_{(n{-}1)k{+}1}^{(k)}$ in the same way as $T_n$ is related to
the order complex of $\Pi_n$, i.e., is $T_n^{(k)}$ homeomorphic to
$\Delta(\Pi_{(n{-}1)k{+}1}^{(k)})$? More than that, can
$\Delta(\Pi_{(n{-}1)k{+}1}^{(k)})$ be obtained from $T_n^{(k)}$ by a
sequence of stellar subdivisions?
\end{qstrm}

An approach to this question along the lines of
Section~\ref{sec_treesAn} does not work right away: The poset
$\Pi_{(n{-}1)k{+}1}^{(k)}$ is not a lattice, and a concept of nested
sets for more general posets is not (yet) at hand.

%%%%%%%%%%%%%%%%%%%%%%%%%%%%%%%%%%%%%%%%%%%%%%%%%%%%%%%%%%%%%%%%%%%%%%%% 
%%%%%%%%%%%%%%%%%%%%%%%%%%%%%%%%%%%%%%%%%%%%%%%%%%%%%%%%%%%%%%%%%%%%%%%% 

\bigskip
The generalization of $T_n$ to complexes of $k$-trees $T_n^{(k)}$ thus
turns out to be somewhat unnatural from the point of view of nested set
constructions. We propose another generalization which is
motivated by starting with a natural generalization of the partition
lattice that remains within the class of lattices.

\begin{dfrm} 
For $n\,{>}\,k\,{\geq}\,2$, the {\em $k$-equal lattice\/} 
$\Pi_{n,k}$ is the
sublattice of the partition lattice $\Pi_n$ that is join-generated by
partitions with a single non-trivial block of size $k$.
\end{dfrm}

Observe that we retrieve $\Pi_n$ for $k{=}2$. There is an extensive
study of the $k$-equal lattice in the literature, mostly motivated by
the fact that $\Pi_{n,k}$ is the intersection lattice of a natural
subspace arrangement, the {\em $k$-equal arrangement\/}.  Its homology has
been calculated by Bj\"orner and Welker~\cite{BWe}, it was shown to be
shellabe by Bj\"orner and Wachs~\cite{BWa2}, and its
$\Sigma_n$-representation theory has been studied by Sundaram and
Wachs~\cite{SWa}.

The irreducibles $\II$ in $\Pi_{n,k}$ are partitions with exactly one
non-trivial block, this time of size at least~$k$. Sets of
irreducibles are nested if and only if for any two elements the
non-trivial blocks are either contained in one another or disjoint.

Constructing trees from nested sets, analogous to the construction of
$T_n$ from $\NN(\Pi_n,\II)$ in the proof of Theorem~\ref{thm_TnminN},
suggests the following definition:

\begin{dfrm}
  The {\em complex of $k$-equal trees\/} $T_{n,k}$ is a simplicial
  complex with maximal simplices  given by combinatorial types of 
  rooted trees $T$ on $n$ labelled
  leaves which are binary except at preleaves, where they are $k$-ary.
  Here, {\em preleaves\/} of $T$ are leaves of the tree that is obtained from
  $T$ by removing the leaves. Lower dimensional simplices are obtained by
  contracting internal edges.
\end{dfrm}

\begin{figure}[ht]
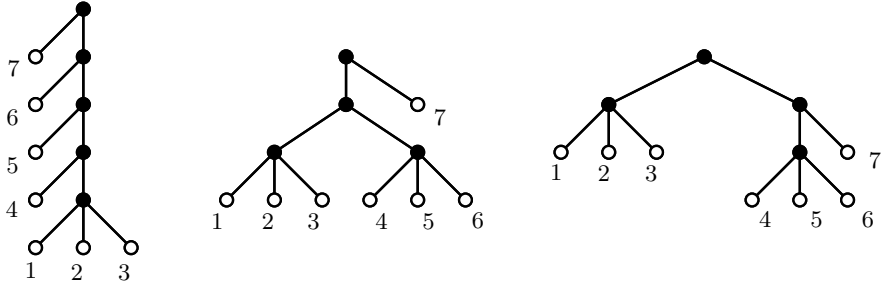

  \begin{picture}(0,0)%
    \includegraphics{kequaltrees.pstex}%
  \end{picture}%
  \input{kequaltrees.pstex_t}%
  
\caption{Maximal simplices in $T_{7,3}$}
\label{fig_kequaltrees}
\end{figure}

We depict the tree types occurring as maximal simplices of $T_{7,3}$
in Figure~\ref{fig_kequaltrees}. Observe that one is a $3$-dimensional
simplex in $T_{7,3}$, whereas the two others are $2$-dimensional. The
definition of $T_{n,k}$ does not appeal as natural, however, it is a
trade off for the following Proposition and Corollary which are
obtained literally for free, having the arguments of
Section~\ref{sec_treesAn} at hand.

\begin{prop}
The complex of $k$-equal trees, $T_{n,k}$,
and the minimal nested set complex of the $k$-equal 
lattice, $\NN(\Pi_{n,k}, \II)$, coincide as abstract simplicial 
complexes. 
In particular, the order complex $\Delta(\Pi_{n,k})$ can be
obtained from $T_{n,k}$ by a sequence of stellar subdivisions.  
\end{prop}

\begin{pf}
  There is a bijection between trees in $T_{n,k}$ and nested sets in
  $\NN(\Pi_{n,k}, \II)$ analogous to the bijection between $T_n$ and
  $\NN(\Pi_n,\II)$ that we described in the proof of Theorem~\ref{thm_TnminN}.
  Referring again to
  \cite[Thm.\ 4.2]{FM}, the complexes are connected by a sequence of 
  stellar subdivisions. 
\end{pf}

\begin{crl} The graded homology groups of the complex of $k$-equal 
trees and of the $k$-equal lattice are isomorphic as 
$\Sigma_n$-modules:
\[
     \wti H_*(T_{n,k}) \,\, \cong_{\Sigma_n} \, \, \wti H_*(\Pi_{n,k})
\]
\end{crl}

%%%%%%%%%%%%%%%%%%%%%%%%%%%%%%%%%%%%%%%%%%%%%%%%%%%%%%%%%%%%%%%%%%%%%%%
%%%%%%%%%%%%%%%%%%%%%%%%%%%%%%%%%%%%%%%%%%%%%%%%%%%%%%%%%%%%%%%%%%%%%%%
%%%%%%%%%%%%%%%%%%%%%%%%%%%%%%%%%%%%%%%%%%%%%%%%%%%%%%%%%%%%%%%%%%%%%%%

%%%%%%%%%%%%%%%%%%%%%%%%%%%%%%%%%%%%%%%%%%%%%%%%%%%%%%%%%%%%%%%%%%%%%%%
%%%%%%%%%%%%%%%%%%%%%%%%%%%%%%%%%%%%%%%%%%%%%%%%%%%%%%%%%%%%%%%%%%%%%%%
%%%%%%%%%%%%%%%%%%%%%%%%%%%%%%%%%%%%%%%%%%%%%%%%%%%%%%%%%%%%%%%%%%%%%%%

\end{document}